\documentclass[graybox]{svmult}

\usepackage{type1cm}        %
\usepackage{makeidx}         %
\usepackage{graphicx}        %
\usepackage{multicol}        %
\usepackage[bottom]{footmisc}%

\usepackage{newtxtext}       %
\usepackage[varvw]{newtxmath}       %
\usepackage{algorithm2e}
\SetKwBlock{Repeat}{repeat}{}
\usepackage{wrapfig}

\usepackage{pgfplots}
\DeclareUnicodeCharacter{2212}{−}
\usepgfplotslibrary{groupplots,dateplot}
\usetikzlibrary{patterns,shapes.arrows}
\pgfplotsset{compat=newest}

\usepackage{hyperref}
\usepackage{background}

\backgroundsetup{
	pages=all,
	color=black,
	placement=top,
	scale=1.0,
	opacity=1,
	vshift=-0.5cm,
	contents={
		\begin{minipage}{1.5\linewidth}
			\scriptsize
			This is a preprint of the following chapter: Benjamin Berkels and Peter Binev, \emph{Joint denoising and line distortion correction for raster-scanned image series}, published in \emph{Multiscale, Nonlinear and Adaptive Approximation II}, edited by Ronald DeVore and Angela Kunoth, 2024, Springer.
            It is the version of the authors' manuscript prior to acceptance for publication and has not undergone editorial and/or peer review on behalf of the Publisher (where applicable).
            The final authenticated version is available online at: \url{https://doi.org/10.1007/978-3-031-75802-7_6}.
		\end{minipage}
	}
}

\makeindex             %

\begin{document}

\title*{Joint denoising and line distortion correction for raster-scanned image series}
\author{Benjamin Berkels and Peter Binev}
\institute{Benjamin Berkels \at RWTH-Aachen, Germany, \email{berkels@aices.rwth-aachen.de}
\and Peter Binev \at University of South Carolina, Columbia, SC 29208, USA, \email{binev@math.sc.edu}}
\maketitle

\abstract{The problem of noise in a general data acquisition procedure can be resolved more accurately if it is based on a model that describes well the distortions of the data including both spatial and intensity changes. The focus of this article is the modeling of the position distortions during sequential data acquisitions. A guiding example is the data obtained by Scanning Transmission Electron Microscopy (STEM) and High Angular Annular Dark Field (HAADF) data, in particular.   
The article discusses different models of the position noise and their numerical implementations comparing some computational results.}

\section{General Data Distortion Model}
\label{sec:1}
The process of recovering information from data usually involves several steps. A standard presumption is to assume that the data set $d$ is a result of applying an operator $\Phi$ to an ideal set $\mathcal{P}$ describing the phenomenon. The goal is to recover $\mathcal{P}$ from $d$ by inverting $\Phi$. Often, the input and output spaces of $\Phi$ are different, e.g., in the case $\Phi$ is a ``parameter-to-state'' map in which the goal is to recover from $d$ the set of parameters $\mathcal{P}$ that determine it.
The mathematical problem of inverting the ``parameter-to-state'' map is well understood in case the data is clean from any distortions typically accompanying the data acquisition process. 
This article considers the case of having an operator $\Phi$ with the same input and output spaces $\Omega$ in which the goal is to recover a ``clean'' version $c$ of the data set $d$. In this case, $\Phi$ is a \emph{distortion} operator acting ``state-to-state''. We consider this as a standalone problem, although anticipate the scenario of it being the first step of the two-step inversion problem recovering $c$ from $d$ which is followed by the ``state-to-parameter'' map recovering $\mathcal{P}$ from $c$.

The choice of the data space $\Omega$ depends on the nature of the distortion process described by the operator $\Phi$. One typical example is the case in which $\Omega$ is the space of images. The continuous version of it is the space $\Omega$ of all bounded functions  $f:[0,1]^2 \rightarrow \mathbb{R}_+^n$, where $\mathbb{R}_+$ is the set of non-negative real numbers. The most studied situation is $n=1$ and we will concentrate on it. 
A straightforward generalization of the methods often works for $n>1$ but the problems become more challenging in some settings for $n>1$, in particular for multimodal data that might even require replacing the range $\mathbb{R}_+^n$ of $f$ with a more intricate set. Other choices of $\Omega$ that can be treated similarly include bounded functions  $f:[0,1]^3 \rightarrow \mathbb{R}_+^n$ representing data sets in 3D space and $f:[0,1]^4 \rightarrow \mathbb{R}_+^n$ for 4D space/time. The only significant problem in such generalizations is the implementation of the algorithms that have to manage the increased amounts of data, memory, and computations.

An abstract problem formulation involves a data function $g:=\Phi\circ f$ that is assumed to be received after the application of the operator $\Phi$ that is unknown or at least has a stochastic component. In reality, we do not have $g$ but some data points that can be treated as the values $\hat g$ of some functionals applied to $g$. Note that $\hat g$ belongs to a finite dimensional space while $\Omega$ is infinite dimensional. On the other hand, $f$ is also approximated by some $\hat f$ that belongs to a chosen class of approximating functions that is finite dimensional. In addition, the operator $\Phi$ is also unknown and has to be approximated by an operator $\hat\Phi$ based on a model depending on finite number of parameters. 
The result of these approximations is a discretized version $\hat g$ of $\hat\Phi\circ\hat f$ that we have to compare to the actual data $d$. 

The problem of interest is typical for sequential data acquisition in which the data is scanned at different positions $p_j$ on a given lattice, assuming that $j\in J$ for a finite index set $J$. The process is time-dependent meaning that both $f$ and $p_j$ change with time. It is often the case that $f$ does not change significantly during scanning and its changes can be attributed to a general 
component.
Then the assumption is  that $f$ does not change during the data acquisition but  instead of sampling $f$ at $p_j$ the data is received by sampling $f$ at distorted positions $\hat p_j$.  

We measure the discrepancy using a distance functional $\operatorname{Dist}_\text{dat}(d,\hat g)$ that is often a square of the $\ell_2$-norm of the difference of the numerical values.
The process of finding  $\hat f$ involves a regularized minimization of this distance. A more complicated distance functional is used in Subsection \ref{sec:2.2.5}.
In the following, we derive it from an assumption of Poisson noise on the input data. Recalling the Poisson distribution $P_\lambda(k)=\frac{\lambda^k}{k!}\exp{(-\lambda)}$, given a true intensity $u\geq0$, the probability to observe the intensity $v\in\mathbb{N}^0$ is given by
\[P(v\vert u)=P_u(v)=\frac{u^v}{v!}\exp{(-u)}\]
Thus, given an observation $v$, the solution $u$ should maximize the probability $P(v\vert u)$.
Therefore, we assign $-\log P(v\vert u)$ as cost for $u$ given $v$, i.e.
\begin{align*}
\operatorname{Dist}_\text{dat}(u,v)&=-\log (P(v\vert u))=-\log \frac{u^v}{v!} +u\\
&=u-v\log u+\underbrace{\log(v!)}_{\text{independent of }u}
\end{align*}
Thus, we get $\operatorname{Dist}_\text{dat}(u,v)=u-v\log u$. Note that this approach can be applied to other noise models as well, as long as the distribution is known.

Both, the operator $\hat\Phi$ and the approximating function $\hat f$, can be forced to belong to appropriate smoothness classes.  
The preferred way to express this enforcement is through two other functionals. 
The first one is motivated by the assumption that the distortion $\hat\Phi$ should be close to the identity mapping and we have a way to measure this closeness. Thus, it is (a power of) a distance $\operatorname{Dist}_\text{op}(\hat\Phi,I)$ between $\hat\Phi$ and the identity $I$ in a chosen metric. The second one is usually (a power of) a norm $\|\hat f\|_Y$ in a chosen Banach space. Often $\|\hat f\|_Y$ is replaced by a semi-norm or a quasi-norm that better expresses the nature of $f$. 
It should be noted that $\operatorname{Dist}_\text{op}(\hat\Phi,I)$ and $\|\hat f\|_Y$ are acting on discrete objects inspired by the corresponding continuous models.

To simplify the notation, we recognize that $g$, $f$, and $\Phi$ are not available to us and will use the same notation for their approximations identifying $\hat g$ with $g$, $\hat f$ with $f$, and $\hat\Phi$ with $\Phi$ from now on. However, we will keep separate notations for $p_j$ and $\hat p_j$.

We find the approximations of $f$ and $\Phi$ using regularized minimization process with a {\em loss function} 
\begin{equation}
\mathcal{L}=\mathcal{L}(\mathcal{P}_f,\mathcal{P}_\Phi)=\operatorname{Dist}_\text{dat}(d,\Phi\circ f) +\lambda\operatorname{Dist}_\text{op}(\Phi,I)+\mu\|f\|_Y  ~,  
\label{loss}
\end{equation}
where $\mathcal{P}_f$ and $\mathcal{P}_\Phi$ are the parameters determining the approximations of $f$ and $\Phi$, while $\lambda\geq0$ and $\mu\geq0$ are tuning parameters weighting the regularization terms for $f$ and $\Phi$.
Sometimes, one or even both of the constants $\lambda$ and $\mu$ are set to zero. This effectively leaves only implicit regularization imposed via the approximation procedure and restrictions on the sets $\mathcal{P}_f$ and/or $\mathcal{P}_\Phi$. 

The simplest interpretation of the above model assumes acquisition of a single frame of data at a rectangular set $P$ of positions $p_j$. Then the minimization of the loss function (\ref{loss}) provides good results only in cases of small distortions and regular data patterns. However, the problems of interest rarely possess either of these features. To process more involved problems one needs to obtain additional information about the distortions. 
One way of gathering this information is to acquire several frames from the same function $f$ assuming that $f$ itself does not change significantly from frame to frame.
This will allow a comparison of the distortions for each individual frame and provide a basis for a better estimation of the approximations of these distortions.
In the next section, we consider different models that take advantage of the multi-frame setup.

\section{Distortion Estimation in the Multi-Frame Setup}
\label{sec:2}
Let $K$ be the number of frames gathered in the experiment. Each frame is a discrete approximation of the ideal image $f$. While the intention is to sample $f$ at the same set $P$ of positions $p_j$ for $j\in J$, in reality $f$ is sampled at different positions $\hat p_j^{[k]}$  for $j\in J$ at different frames $k=1,2,...,K$. 
Each frame is determined by the set of distorted positions $P^{[k]}:=\{\hat p_j^{[k]}, j\in J\}$ corresponding to the discretization $g_k:=f(P^{[k]})$.

The general model must be slightly modified to account for the multiple frames. The operator $\Phi$ is now acting from $\Omega$ to $\Omega^K$ and the regularization through $\operatorname{Dist}_\text{op}$ should compare $\Phi$ with $I^K$ instead of $I$. The data $d=\{d_k\}_{k=1}^K$ is also treated as $K$-dimensional entity.
In addition, instead of assuming that $\Phi$ is such that $\Phi\circ f=g$ as in the general model, we will focus on its position distortion part 
and will assume that $\Phi\circ f=\left(g_k\right)_{k=1}^K$ anticipating that the estimation of the other distortions will be handled by the regularization term $\operatorname{Dist}_\text{dat}(d,\Phi\circ f)$ by penalizing the distances between $d_k$ and $g_k$.
We associate the operator  $\Phi$ with a vector of operators $\left(\phi_k\right)_{k=1}^K$ that are transforming $P^{[k]}$ back to $P$ yielding  $g_k\circ\phi_k = f$.
In some sense $\left(\phi_k\right)_{k=1}^K$ represents the inverse operator of $\Phi$.
To state the updated form of (\ref{loss}), we need to redefine the two distance functionals for $\Phi$ as sums of individual distance functionals $\operatorname{Dist}_\text{op}^\diamond$ and $\operatorname{Dist}_\text{dat}^\diamond$ for each frame and $\phi_k$ reminding the reader that often for the formulation of a loss function the distance functionals are powers of an actual distance measure.  
We set
\[\operatorname{Dist}_\text{op}(\Phi,I^K):=\sum_{k=1}^K \operatorname{Dist}_\text{op}^\diamond (\phi_k,I)
\]
and 
\[\operatorname{Dist}_\text{dat}(d,\Phi\circ f):=\sum_{k=1}^K \operatorname{Dist}_\text{dat}^\diamond (d_k\circ\phi_k, f).
\]
The general form of a loss function for the multi-frame setup is
\begin{equation}
\mathcal{L}=\mathcal{L}(\mathcal{P}_f,\mathcal{P}_\Phi)=\sum_{k=1}^K \left(\operatorname{Dist}_\text{dat}^\diamond (d_k\circ\phi_k, f) +\lambda \operatorname{Dist}_\text{op}^\diamond (\phi_k,I)\right)+\mu\|f\|_Y .  
\label{loss2}
\end{equation}

Depending on the settings, the minimization of the loss function could be a very complicated process. Therefore, simplifying the dependencies in the formulation could be beneficial. In general, $\operatorname{Dist}_\text{op} (\Phi_k,I)$ in (\ref{loss}) could depend on $f$ since $\Phi: \Omega \to \Omega$ is defined for $\Omega$ being a set of functions. However, in (\ref{loss2}) we define $\phi_k$ as a transformation between sets of positions avoiding a potential dependence on $f$. 

There are several possible applications of the models described by the loss function (\ref{loss2}). One of the most interesting is scanning transmission electron microscopy (STEM). In the next subsection, we consider the specifics in the STEM applications.

\subsection{Atomic Scale Data Processing in Scanning Transmission Electron Microscopy (STEM)}
\label{sec:2.1}
STEM images are usually two-dimensional rectangular arrays of counts of the numbers of electrons that have hit the detector(s) during the short period of time when the electron beam has been positioned at the corresponding spot. The interactions of the electrons from the beam with the atoms of the specimen result in changing the trajectories of some electrons. The standard detectors are \textit{bright field} (BF) detecting the electrons whose trajectories are not changed substantially, \textit{annular bright field} (ABF) detecting the electrons diverted by a small angle, \textit{annular dark field} (ADF) detecting the electrons diverted by a moderate angle, and \textit{high angle annular dark field} (HAADF) detecting the electrons diverted by a high angle. HAADF is the most used STEM acquisition mode due to the so-called Z-contrast allowing to better distinguish atoms with different number of protons. The examples in this article are for HAADF data at atomic scale resolution. As discussed in the introduction, one can also consider a vector data output ($n>1$), e.g., by combining the data of BF, ABF, ADF, and/or HAADF or even consider a two-dimensional array of detectors resulting in so-called 4D-STEM data (see, e.g. \cite{Ophus_2019}). While the vector data would result in better estimation of the distortions, the associated algorithms are more complicated and their computational complexity is higher.

The electron beam is interacting with the specimen within a small volume around the central axis of the beam.
The density of the interaction is higher closer to the axis. It is usually modeled by a 2D Gaussian distribution on the planes perpendicular to the axis although in actuality it is much more complicated.
We identify the axis with its intersection point $p$ with the $xy$-plane of the specimen.
The data collected for the corresponding volume of interaction is the (approximate) value of $f(p)$.
The standard scanning routine moves the electron beam in $x$-direction for a given number of equal steps $h_x$ at fixed $y$-position and then changes the $y$-position by a step $h_y$ and starts scanning from the initial $x$-position. The scanning continues for a given number of steps in $y$-direction. 
The electron beam is positioned by incremental changes of an electromagnetic field guiding it. Different environmental factors can influence this electromagnetic field and distort the intended position of the beam - the higher the resolution, the bigger the relative distortions. 
In atomic scale resolution, these position distortions are the main obstacle of receiving higher precision. 
This is the main motivation for developing acquisition procedures and numerical algorithms for estimating the position distortions.
The beam is shut out during the repositioning but might not be completely steady in maintaining the same position during the acquisition process of a single data point (or data vector). 
This would change somewhat the density of interaction within the volume corresponding to $p$ creating data that is deviating from $f(p)$. 
Such a deviation is expected to be more pronounced for the first point with a new $y$-position due to the larger repositioning from the previously scanned point. 
As our goal is the estimation of the position distortions, this and other value distortions will be handled by the regularization terms of the loss function.

\subsection{Models and algorithms for estimating the position distortion in HAADF-STEM data}
\label{sec:2.2}
We consider five particular models. The first two are based on continuous models that are discretized and the others are discrete by design.
In Subsection \ref{sec:2.2.5} we introduce a new discrete model for joint denoising and distortion correction.

\subsubsection{Iterative approximating algorithm based on the estimation of the position distortions of pairs of consecutive frames}
\label{sec:2.2.1}
The method introduced in \cite{BeBiBl13,YaBeDa14} does not consider a global minimization procedure for a given loss function like (\ref{loss2}) but suggests an iterative process consisting of several steps.
The input is the data for the frames $d_k$, $k=1,...,K$ and the output is the approximation $f$ calculated  together with position distortions  $\phi_k$, $k=1,...,K$. 
The key ingredient of the method is the minimization process that for given two functions $u$ and $v$ defined on a domain $\mathcal{D}$ finds a transformation $\phi$ for which the functions $u(p)$ and $v\circ \phi(p):=v(\phi(p))$ are as close as possible in terms of minimizing the loss function
\begin{equation}
\mathcal{L}^\star=\mathcal{L}^\star(\mathcal{P}_\phi)= \operatorname{Dist}_\text{dat} (u,v\circ \phi) +\lambda R(\phi)  
\label{lossd}
\end{equation}
among all vector-valued piecewise bilinear functions $\phi : \mathcal{D} \to \mathcal{D}$ with a positive semi-definite Jacobian $\mathcal{J}(\phi)$. 
The similarity measure $\operatorname{Dist}_\text{dat}^\diamond$ is based on the normalized cross correlation and is defined by
\begin{equation} 
\operatorname{Dist}_\text{dat}^\diamond(u,v):= - \frac{1}{|\mathcal{D}|}\int_\mathcal{D} \frac{u-\Bar{u}}{\sigma_u} ~ \frac{v-\Bar{v}}{\sigma_v} ~\mathrm{d}p ~,
\label{dist}
\end{equation}
where as usual $\Bar{u}:=(1/|\mathcal{D}|)\int_\mathcal{D} u(p) ~\mathrm{d}p$ is the mean value of $u$ over the domain $\mathcal{D}$ and 
$\sigma_u:=\sqrt{(1/|\mathcal{D}|)\int_\mathcal{D} |u(p)-\Bar{u}|^2~\mathrm{d}p}$ is the standard deviation of $u$.
The regularization term $R(\phi)$ is the squared integrated Frobenius norm of the difference between the Jacobian $\mathcal{J}(\phi)$ of $\phi$ and the identity $I$: 
\begin{equation} 
R(\phi):=\int_\mathcal{D}\|\mathcal{J}(\phi)-I\|_F^2 \mathrm{d}p. 
\label{R}
\end{equation}
It is corresponding to the term $\operatorname{Dist}_\text{op}^\diamond(\phi_k,I)$ in (\ref{loss2}).
The parameter $\lambda$ is used as a normalization that makes the two terms in (\ref{lossd}) of similar magnitude.

In the above description, it is assumed that both functions $u$ and $v$ are defined in the same domain $\mathcal{D}$. If this is not the case and we have that $u$ is defined in $\mathcal{D}_u$ and $v$ is defined in $\mathcal{D}_v$, then we can use $\mathcal{D}:=\mathcal{D}_u\cap\mathcal{D}_v$. However, we can opt for using $\mathcal{D}:=\mathcal{D}_u\cup\mathcal{D}_v$ by first calculating the mean values $\Bar{u}:=(1/|\mathcal{D}_u|)\int_{\mathcal{D}_u} u(p) ~\mathrm{d}p$ and $\Bar{v}:=(1/|\mathcal{D}_v|)\int_{\mathcal{D}_v} v(p) ~\mathrm{d}p$ and then extending $u$ and $v$ by their mean values at the points they are not defined.
Note that the product under the integral in (\ref{dist}) is zero for $p\notin \mathcal{D}_u\cap\mathcal{D}_v$.

A key feature of the algorithm is using that the distortion between two consecutive frames is close to the identity and it can be used as initial approximation. 
The loss function (\ref{lossd}) is defined for continuous functions and thus we need to define functional analogs $g_k$ of the discrete data $d_k$. It can be done in different ways but the most convenient one is to define $g_k$ to be the piecewise bilinear function on a rectangular grid interpolating the data $d_k$ at the integer points. 
As a preprocessing step of the iterative procedure, we determine the transformations $\phi=\phi_{k+1,k}$ that make the functions
$u=g_k$ and $v\circ\phi=g_{k+1}\circ \phi_{k+1,k}$ as close as possible by minimizing (\ref{lossd}) with initial approximation $\phi=I$.
The practical implementation of the minimization is rather involved. One of the issues that should be addressed is the near-periodicity of the data that creates several local minima. This necessitates involvement of multiresolution processing and estimation of the rigid registration component in the procedure. Details can be found in \cite{BeBiBl13,YaBeDa14}. 

The outline of the general algorithm is as follows:

\begin{algorithm}[H]
    \TitleOfAlgo{Iterative approximation using general bilinear distortions}
	\KwIn{$g_k$ for $k=1,...,K$ and $\phi_{k+1,k}$ for $k=1,...,K-1$;}
	\textbf{Initialization:}  set the initial guess $f^{[0]}=g_1$ and $\phi_1=I$\;
    \Repeat{estimate $\phi_1$ by the approximate minimizer $\phi$ of (\ref{lossd}) with $u=f^{[0]}$ and $v=g_1$ using as initial guess the last version of $\phi_1$\;
    consecutively for $k=2,...,K$ estimate $\phi_k$ by the approximate minimizer $\phi$ of (\ref{lossd}) with $u=f^{[0]}$ and $v=g_k$ using as initial guess $\phi_{k,k-1}\circ\phi_{k-1}$\;
    calculate $f$ as the average (e.g., median or mean) of the functions $g_k \circ \phi_k$ for $k=1,...,K$\;
    \eIf{the error estimator $\|f-f^{[0]}\|$ is below a prescribed tolerance or other stopping criterion is satisfied (number of iterations, time of computation, etc.). \; }{break\;}
    {set $f^{[0]}:=f$\;}
    }
\end{algorithm} 

There can be different variants of this algorithm depending on the type of average used and the way the domains $\mathcal{D}_k$ of each of the frames are treated. 
Usually, the median averaging produces better results. 
The standard approach is to consider the intersection of the domains $\mathcal{D}:=\cap_{k=1}^K\mathcal{D}_k$ but it is often beneficial to work with weighted norms over the union $\cup_{k=1}^K\mathcal{D}_k$ assigning higher weights on the subregions belonging to more domains $\mathcal{D}_k$.  
  
\subsubsection{Iterative approximating algorithm with global regularization}
The iterative algorithm proposed in Subsection \ref{sec:2.2.1} attempts to approximate well the solution of the global minimization problem based on (\ref{loss2}). 
It performs the minimization of (\ref{lossd}) frame by frame for a fixed function $f=f^{[0]}$ that is updated at each iteration by an average of $g_k \circ \phi_k$ for $k=1,...,K$. 
The averaging process can be adjusted to yield a function $f$ with specific properties which are usually imposed by the term $\mu\|f\|_Y$ in (\ref{loss2}).
Using weighted averaging can balance well the terms $\operatorname{Dist}_\text{dat}^\diamond(f,g_k \circ \phi_k)$ in the corresponding sum in (\ref{loss2}) making it close to optimal.
However, the sum $\sum_{k=1}^K \lambda \operatorname{Dist}_\text{op}^\diamond (\phi_k,I)= \lambda \sum_{k=1}^K R(\phi_k)$ might stay far from optimal for significant number of iterations.

The algorithm proposed in \cite{BeLi19} adds a step at the end  of each iteration that improves substantially the optimality of $\sum_{k=1}^K R(\phi_k)$ leading to very good results after a few iterations.
The idea is to introduce a local change of variables $\psi$ replacing $f$ with $f\circ\psi$ and $\phi$ with $\phi\circ\psi$ in the sums of the loss function.
The values of terms $\operatorname{Dist}_\text{dat}^\diamond(f\circ\psi,g_k \circ \phi_k\circ\psi)$ are only slightly perturbed since they are calculated by (\ref{dist}) with integration over $\psi(\mathcal{D})$ instead of $\mathcal{D}$ and $\psi(p)\approx p$.
The $\psi$ is calculated by minimizing $\sum_{k=1}^K ||\phi_k\circ\psi-\operatorname{id}_{\mathcal{D}}||_{L^2}^2$ for the estimated distortions $\phi_k$ at the current iteration, where $\operatorname{id}_{\mathcal{D}}$ is the identity mapping on $\mathcal{D}$.
Then $f$ is updated as an average of the functions $g_k \circ \phi_k\circ\psi$ for $k=1,...,K$.
As already mentioned, this modification of the algorithm from Subsection \ref{sec:2.2.1} converges very fast to the desired solution and noticeably reduces the bias of the reconstruction $f$ towards the coordinate system of $g_1$, which is due to the alternating minimization scheme used in \cite{BeBiBl13} combined with using $g_1$ as initial guess for $f$.
\subsubsection{Smart align}
In \cite{JoYaPe15}, the authors suggest a non-rigid registration methodology and a software package called Smart Align. The procedure uses two types of images, a static reference image and a set of ``moving'' images.
The reference image is calculated initially as an average of the rigidly aligned frames. Once the moving images are updated, it is recalculated as an average of them to start the next iteration of the procedure. 
The original frames are the starting state of the moving images each of which is gradually modified by a gradient descent process utilizing the difference of the gradients of a moving image and the reference image to detect the position varying offsets between them. This results in obtaining an incremental calculation field
that is improved through implementing constrains motivated by the specifics of the STEM images.
The software presents solutions to some practical situations like processing frames with beam damage. 
Smart Align presents an alternative approach to the loss function minimization discussed in the previous subsections providing a good algorithmic solution to the distortion correction problems.

\subsubsection{Joint denoising and distortion correction with a bump fit image model}
\begin{figure}[t]
\centering
\begin{tabular}{@{}cc@{}}
\includegraphics[width=0.45\linewidth]{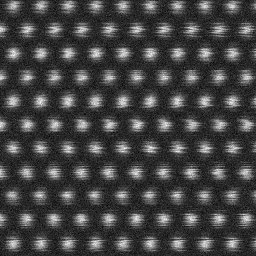}&
\includegraphics[width=0.45\linewidth]{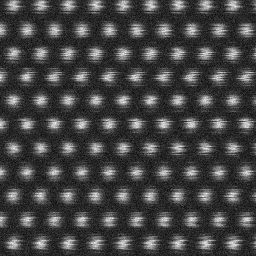}\\
$d_1$&$d_{64}$
\end{tabular}
\caption{\label{fig:SyntheticInput}First and last image of the synthetic STEM images series used for the numerical performance evaluation. %
}
\bigskip
\begin{tabular}{@{}cc@{}}
\includegraphics[width=0.45\linewidth]{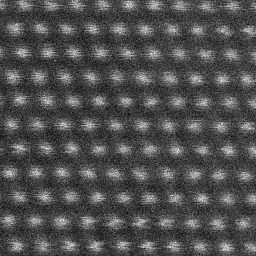}&
\includegraphics[width=0.45\linewidth]{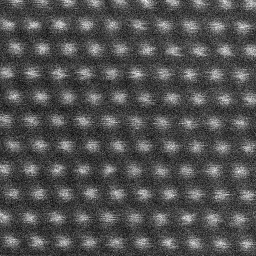}\\
$d_1$&$d_{64}$
\end{tabular}
\caption{\label{fig:RealInput}First and last image of the subset of the HAADF STEM GaN $[11\overline{2}0]$ image series from \cite{YaBeDa14} used for the numerical performance evaluation here.
}
\end{figure}

In \cite{BeWi17}, a Bayesian method is derived that, for a given STEM image series, jointly estimates the distortion in each image and reconstructs the underlying image of the material.
For this discrete model, the data $d=\{d_k\}_{k=1}^K$ is assumed to be given as $K$ noisy STEM pixel images of size $M\times N$ each, i.e., $d_k\in\mathbb{R}^{M\times N}$, $k=1,\ldots,K$, cf. Figure~\ref{fig:SyntheticInput}. Here, the pixel values are considered to correspond to a Cartesian pixel grid discretizing the image domain $\mathcal{D}=[0,1]^2\subset\mathbb{R}^2$. The coordinates of the pixels are denoted by $x_{ij}\in\mathcal{D}$, $(i,j)\in\{1,\ldots,M\}\times\{1,\ldots,N\}$.

The sought denoised, continuous image $f:\mathcal{D}\to\mathbb{R}$ is modelled as a sum of parametrized Gaussians, i.e.
\begin{equation*}
f[\mathcal{P}](x)=\sum_{l=1}^L b[c_l](x-y_l)+o\,,
\end{equation*}
where $L$ denotes the number of atoms, $b$ is a parameterized Gaussian, $c_l$ and  $y_l$ are parameters and the center of the Gaussian for atom $l$, respectively, and $o$ is the background intensity. $\mathcal{P}=(L,y_1,\ldots,y_L,c_1,\ldots,c_L,o)$ is the vector of all parameters.

The deformation aligning $f$ with $d_k$ is split into a global sample drift in terms of a rigid body motion consisting of a rotation $R_k\in SO(2)$ and a translation $v_k\in\mathbb{R}^2$ (using the positions of $d_1$ as reference coordinate system) and the local displacement/shift $s_{ij}^k\in\mathbb{R}^2$ for each pixel $i,j$. Note that this shift is linked to a scanline and thus it is inherently discrete in $y$-direction.
Thus, the sought image $f$ approximates the discrete image $d_k$ via
\begin{equation*}
f[\mathcal{P}](R_k(x_{ij}+s_{ij}^k)+v_k)\approx (d_k)_{ij}\,.
\end{equation*}
The objective functional is
\begin{align*}
\mathcal{L}^K[&\mathbf{s}^1,\ldots,\mathbf{s}^K,\mathcal{P},(R_k,v_k)_{k=1}^K]
\\&=\sum_{k=1}^K\bigg[
\sum_{i=1}^M\sum_{j=1}^N \operatorname{Dist}_\text{dat}\left(f[\mathcal{P}](R_k(x_{ij}+s_{ij}^k)+v_k),(d_k)_{ij}\right)+R_1[\mathbf{s}^k]+R_2[\mathbf{s}^k]
\bigg],%
\end{align*}
where with a specific scaling factor $\Lambda$, $\Delta t$ being the time to scan one pixel, and $\Delta T$ being the time between scanning two consecutive lines 
\begin{equation*}
R_1[\mathbf{s}]=\frac{\Lambda}{2}\left[\frac{|s_{11}|^2}{\Delta t}+\sum_{j=2}^N\left(\frac{|s_{1j}-s_{M(j-1)}|^2}{\Delta T}+\sum_{i=2}^M\frac{|s_{ij}-s_{(i-1)j}|^2}{\Delta t}\right)\right]
\end{equation*}
for $\mathbf{s}\in(\mathbb{R}^2)^{M\times N}$, which results from the assumption that sample undergoes a Brownian motion \cite{BeWi17}. For numerical reasons,
\begin{equation*}
R_2[\mathbf{s}]=\frac1{2}\sum_{i=1}^M\sum_{j=1}^N\left(\nu_\text{hor}|(s_{ij})_1|^2+\nu_\text{vert}|(s_{ij})_2|^2\right)
\end{equation*}
is further added. Here, $(s_{ij})_1$ and $(s_{ij})_2$ denote the horizontal and vertical component of $s_{ij}\in\mathbb{R}^2$ and $\nu_\text{hor},\nu_\text{vert}>0$ are scaling factors.
We would like to point out that the input data $(d_k)_{ij}$ is not interpolated at all. Due to the nature of the STEM acquisition process, interpolation in vertical direction, i.e. between different scanlines, must be avoided when trying to estimate the discontinuous scanline shifts.

The numerical optimization of this model is also rather involved. In particular, getting an initial guess for $\mathcal{P}$ that contains all atoms visible in the reconstructed image is challenging. For details on this, we refer to \cite{BeWi17}. Further down, we describe the other aspects of the minimization strategy, since they will be used for the new model we propose here.

\subsubsection{A new model - Joint denoising and distortion correction with a spline based image model}
\label{sec:2.2.5}

In order to generalize \cite{BeWi17} to images that cannot be described as sum of Gaussians, we propose to use a new model for the sought denoised continuous image $f:\mathcal{D}\to\mathbb{R}$. It is parametrized using the Cartesian product space of 1D cubic B-splines. More specifically, let $\varphi_1,\ldots,\varphi_n$ be a basis of the cubic B-splines on $[0,1]$ with equidistant knots \cite{Bo78}. In other words, the value of $f$ at a position $x=(x_1,x_2)\in\mathcal{D}$ is given by
\[f[\mathcal{P}](x_1,x_2)=\sum_{k,l=1}^np_{kl}\varphi_k(x_1)\varphi_l(x_2),\]
where $\mathcal{P}=(p_{kl})_{k,l=1}^n$ is the parameter vector.

The resulting objective function is numerically minimized using a Trust Region solver \cite{NoWr06}. Due to the non-convex nature of this minimization problem, a careful minimization strategy is necessary. To this end, we follow a similar strategy to the one described in \cite{BeWi17}, which considers a similar approach where $u$ is not represented using splines, but as a sum of Gaussians, one Gaussian for each atom visible in the image domain. $\mathbf{s}^1,\ldots,\mathbf{s}^K,\mathcal{P}$ are each initialized with 0.
The global rotations and translations $(R_k,v_k)_{k=1}^K$ are initialized using the strategy of \cite{BeBiBl13}, but just with a rigid deformation model instead of a non-rigid one. 
Then, $\mathcal{P}$ is updated, by fitting $f[\mathcal{P}]$ to the average of the rigidly aligned input images, i.e. by minimizing
\[\sum_{i=1}^M\sum_{j=1}^N \operatorname{Dist}_\text{dat}\left(f[\mathcal{P}](x_{ij}),\frac1K\sum_{k=1}^Ng_k\big(R_k^{-1}(x_{ij})+v_k^{-1}\big)\right),\]
with respect to $\mathcal{P}$ using the BFGS Quasi Newton algorithm \cite{NoWr06}.
Here, $R_k^{-1}$ and $v_k^{-1}$ denote the rotation and translation corresponding to the inverse of the rigid transformation given by $R_k$ and $v_k$, $g_k$ is the piecewise bilinear interpolation of $d_k$.
With this rough but reasonable guess for $f[\mathcal{P}]$, the objective function is minimized with respect to $\mathcal{P}$ and $(R_k,v_k)_{k=1}^K$ using Trust Region. From here on, the global rigid alignment parameters $(R_k,v_k)_{k=1}^K$ are considered final and thus fixed. 
Then, $\mathcal{L}^K[\mathbf{s}^1,\ldots,\mathbf{s}^K,\mathcal{P}]$ is minimized, while the shifts $\mathbf{s}^k$ are constrained to be constant in each scan line, i.e. $s_{ij}^k=s_{lj}^k$ for all $i,l\in\{1,\ldots,M\}$, $j\in\{1,\ldots,N\}$ using Trust Region. Finally, $\mathcal{L}^K[\mathbf{s}^1,\ldots,\mathbf{s}^K,\mathcal{P}]$ is minimized using Trust Region.

Note that the initialization of this model is simpler than the one for \cite{BeWi17}. For the latter, one needs to know the number of atoms and a good initial guess of the position of each atom. In \cite{BeWi17}, this is obtained by using a special convex surrogate model that treats the detection of the atom as a deconvolution problem with a known kernel (a Gaussian representing a single atom / atomic column).

\section{Numerical results}
\label{sec:NumericalResults}
For the numerical evaluation of the performance of the proposed approach, we created a series consisting of 64 synthetic STEM images. Here the image is modelled as a sum of Gaussians located on a regular grid, intensity noise is added using a Poisson noise model, location noise is added by assuming that the sample is undergoing a Brownian motion. The same approach is used for the generation of synthetic test images in \cite{BeWi17}. Figure~\ref{fig:SyntheticInput} shows the first and the last image from this series. This synthetic data set was modelled to resemble the characteristics of the HAADF STEM GaN $[11\overline{2}0]$ image series used in \cite{YaBeDa14}, cf. Figure~\ref{fig:RealInput}.

To this synthetic dataset, the non-rigid reconstruction \cite{BeBiBl13,YaBeDa14} (denoted by ``NRR''), the extended non-rigid reconstruction with bias correction \cite{BeWi17} (denoted by ``NRR+'') and the method proposed here (denoted by ``JUD'', run with
$\nu_\text{hor} = 25.9$, %
$\nu_\text{vert} = 71.4$, %
$\frac{\Lambda}{2\Delta t} = 50.0$, %
$\Delta T=1000 \Delta t$, $n=2^6$) are applied.

\begin{wrapfigure}[14]{r}{.4\linewidth}
\vspace{-1.5\baselineskip}
\resizebox{\linewidth}{!}{%
\includegraphics{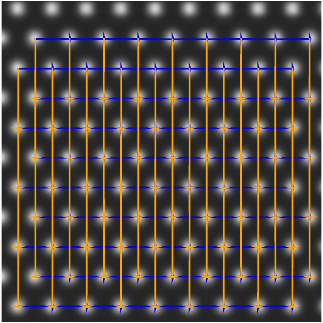}
}%
\caption{\label{fig:SyntheticPrecMatch}Distances used to compute the precision in $x$ and $y$-direction.
}%
\end{wrapfigure}
\begin{figure}[t]
\includegraphics[scale=0.99]{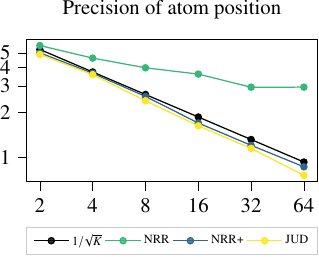}
\includegraphics[scale=0.99]{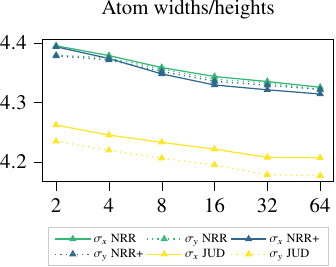}
\caption{\label{fig:SyntheticPrecision}Left: Precision plot for the synthetic test data set for different methods. The horizontal axis shows the number of images used for the reconstruction, the vertical axis the precision value in picometer. Right: Estimated atom sizes in pixels for a synthetic test data set.}
\end{figure}

The quality of the reconstruction here is measured in terms of the so-called \emph{precision}. In electron microscopy, this refers to the precision with which atoms can be located and is computed on images assumed to show a perfect crystal. To compute it, all atoms in a given image are individually fitted using a parameterized Gaussian bump. Then all pairwise distances of neighboring atoms in $x$ and $y$-direction are computed. Finally, the standard deviation of these distances are computed, resulting in a precision in $x$-direction and a precision in $y$-direction. The overall precision is defined as the square root of the sum of the squares of the precision in $x$ and $y$-direction. Since the image is supposed to show a perfect crystal, all pairwise distances in one direction should be the same. In this sense, the precision implicitly measures how well atom centers can be located. 
Figure~\ref{fig:SyntheticPrecMatch} illustrates the pairwise distances used for the precision. $x$-neighboring atoms are connected with blue arrows, $y$-neighboring atoms with yellow arrows. To convert the precision from pixels to picometer, the size a pixel in picometer is estimated as follows. The average of all $x$-distances and all $y$-distances in pixels is computed from all reconstructions. Combined with the expected physical distances for GaN $[11\overline{2}0]$ in $x$ and $y$-directions of
276.174\,pm
and
518.5\,pm%
, respectively, this allows to compute a pixel size in
pm
in $x$ and $y$-direction. Their average is finally used as conversion factor. A very similar conversion logic was used in \cite{YaBeDa14}.
The left plot in Figure~\ref{fig:SyntheticPrecision} shows the resulting precision for the different methods depending on the number of input images.
To compute the precision for a given number of images, the image series is split into non-overlapping series to generate multiple samples. For instance, for two images, the series of 64 images is split into 32 series with 2 images each, each method is applied to each of the 32 small series and the resulting precision values are averaged to get the precision for two images. The only exemption here is the case of one image. Here only every other image is used, i.e., the methods are applied to 32 series with one image each.
As to be expected, the precision of all methods improves with the number of images. However, the scaling is different. Both ``NRR+'' and the proposed method scale as one over the square root of the number of images, ``NRR'' scales noticeable worse. While the proposed method shows a slightly better precision than ``NRR+'', both methods are equally suited when it just comes to locating the atom centers. Looking beyond just the precision, this picture changes though.
The right plot in Figure~\ref{fig:SyntheticPrecision} shows the average sizes of the Gaussians used to compute the precision. Note that the size of the Gaussian used to create the test data is 4.25 pixels. Thus, the proposed method more accurately reconstructs this dataset than the other methods.

\begin{figure}[p]
\centering
\begin{tabular}{@{}cc@{}}
Ground truth&$f$ (NRR)\\
\includegraphics[width=0.44\linewidth]{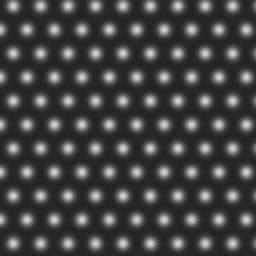}&
\includegraphics[width=0.44\linewidth]{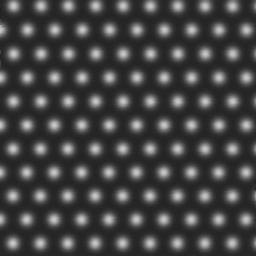}\\
\includegraphics[width=0.44\linewidth]{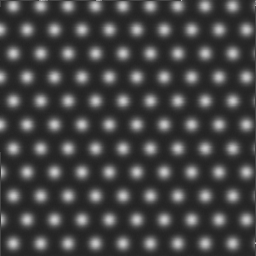}&
\includegraphics[width=0.44\linewidth]{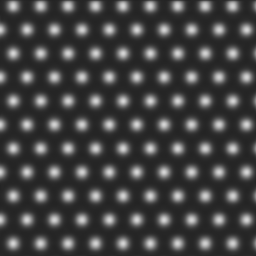}\\
$f$ (NRR+)&$f$ (JUD)\\
\end{tabular}
\caption{\label{fig:SyntheticReconstructions}Ground truth image the synthetic STEM image series was created from and the three reconstructions obtained with the different methods using the whole series of 64 images.}
\medskip
\begin{tabular}{@{}cc@{}}
\includegraphics[width=0.44\linewidth]{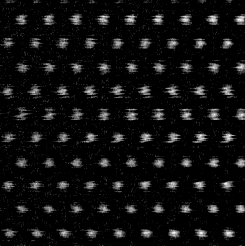}&
\includegraphics[width=0.44\linewidth]{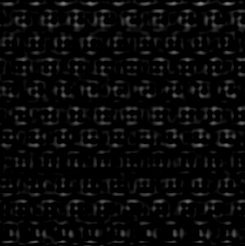}\\
$|f$ (NRR+) - Ground truth$|/|$Ground truth$|$&$|f$ (JUD)- Ground truth$|/|$Ground truth$|$
\end{tabular}
\caption{\label{fig:SyntheticDiff}Pixel-wise relative error of the ground truth image to the reconstructed images from the bottom row of Figure~\ref{fig:SyntheticReconstructions}. $0$~\includegraphics[width=.1\linewidth,height=.5\baselineskip]{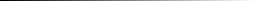}~$5.5\%$}
\end{figure}
Figure~\ref{fig:SyntheticReconstructions} shows the reconstructed images when the whole series of 64 images is used. Side by side, the images look very similar, even though the corresponding precision varies strongly between ``NRR'' on the one hand, and ``NRR+'' and the proposed method on the other hand, cf. Figure~\ref{fig:SyntheticPrecision}, left. Even though the latter two are rather similar in terms of precision (with the proposed method still outperforming ``NRR+''), the proposed method reconstructs the intensity values of the ground truth considerably more accurately than ``NRR+''. This becomes visible when comparing the intensity of the ground truth to the reconstruction, cf. Figure~\ref{fig:SyntheticDiff}. Note that some rows and columns at the boundaries have been cropped to avoid boundary artifacts. It is apparent that the intensity errors of ``NRR+'' are not only larger, but that ``NRR+'' has not fully corrected the scanline artifacts, which is not surprising, since this method does not allow for discontinuous displacements and thus cannot fully correct those.

\begin{wrapfigure}[12]{r}{.5\linewidth}
\vspace{-1.5\baselineskip}
\resizebox{\linewidth}{!}{%
\includegraphics{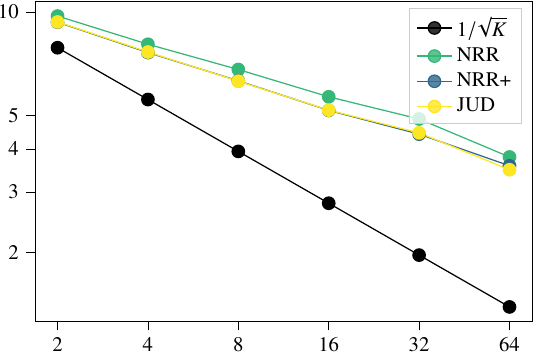}%
}
\caption{Precision plot for the HAADF STEM GaN $[11\overline{2}0]$ image series.}
\label{fig:RealPrecision}
\end{wrapfigure}

\begin{figure}[h]
\centering
\begin{tabular}{@{}ccc@{}}
\includegraphics[width=0.32\linewidth]{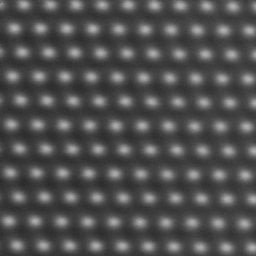}&
\includegraphics[width=0.32\linewidth]{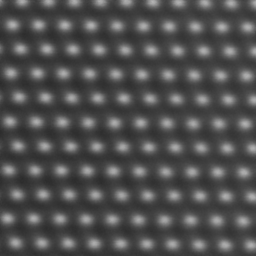}&
\includegraphics[width=0.32\linewidth]{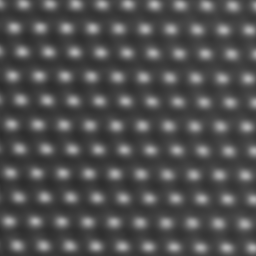}\\
$f$ (NRR)&$f$ (NRR+)&$f$ (JUD)
\end{tabular}
\caption{\label{fig:RealReconstructions}The three reconstructions obtained with the different methods using the whole 64 HAADF STEM GaN $[11\overline{2}0]$ image series.}
\end{figure}

We also applied the three methods on a series of 64 real STEM images. This are the first 64 images of the 512 HAADF STEM GaN $[11\overline{2}0]$ image series used in \cite{YaBeDa14}. Figure~\ref{fig:RealInput} shows the first and the last of these 64 images.
Figure~\ref{fig:RealPrecision} shows the resulting precision in picometer, where the same unit conversion logic as for the synthetic data was used with the same physical reference distances. The relative performance of the different methods is similar to the synthetic data case. The differences however are much smaller and the precision decreases slower than one over the square root of the number of images, which indicates that all of these models still miss some relevant effects encountered in real STEM data.
Note that the precision numbers here are not comparable to the ones listed in \cite{YaBeDa14}, since the reconstructed image for a given number of images were computed very differently. In \cite{YaBeDa14}, series of 512 images were considered and an entire series was aligned and averaged following the algorithm from \cite{BeBiBl13}. Then, the first $n$ of the resulting aligned images were averaged and used to compute the precision value for $n$ frames. Figure~\ref{fig:RealReconstructions} illustrates the reconstructions using all images for the different methods. Like in the synthetic case (cf. Figure~\ref{fig:SyntheticReconstructions}), the reconstructions look very similar. However, since no ground truth is available, we cannot illustrate the effect of the improved scanline distortion correction for the real data case as in the synthetic data case.

\begin{acknowledgement}
The work of B. Berkels was supported by the German research foundation (DFG) within the Collaborative Research Centre SFB 1394 ``Structural and Chemical Atomic Complexity—From Defect Phase Diagrams to Materials Properties'' (Project ID 409476157) in Project A04.
The work of P. Binev was supported by the NSF grants DMS 2038080 and DMS 2245097.
\end{acknowledgement}

\bibliography{references}

\end{document}